\documentclass[12pt]{article}                                               
\input{amssym.def}                                                           
\input{amssym.tex}                                                           
                                                                             
\begin{document}                                                             
%**********************************************                              
\title{Remark on Weil's conjectures}

\author{Igor  ~V. ~Nikolaev}

%**************************************************

\date{}
 \maketitle

%**************************************************

\newtheorem{thm}{Theorem}
\newtheorem{lem}{Lemma}
\newtheorem{dfn}{Definition}
\newtheorem{rmk}{Remark}
\newtheorem{cor}{Corollary}
\newtheorem{exe}{Exercise}
\newtheorem{exm}{Example}

%*************************************************
\newcommand{\Qcoh}{\hbox{\bf Qcoh}}
\newcommand{\QGr}{\hbox{\bf QGr}}
\newcommand{\ch}{\hbox{\bf Char}}
\newcommand{\jac}{\hbox{\bf Jac}}
\newcommand{\ka}{\hbox{\bf k}}
\newcommand{\n}{\hbox{\bf n}}
\newcommand{\mod}{\hbox{\bf mod}}
\newcommand{\Z}{\hbox{\bf Z}}
\newcommand{\Q}{\hbox{\bf Q}}

\newcommand{\Coh}{\hbox{\bf Coh}}
\newcommand{\Mod}{\hbox{\bf Mod}}
\newcommand{\Irred}{\hbox{\bf Irred}}
\newcommand{\Spec}{\hbox{\bf Spec}}
\newcommand{\Tors}{\hbox{\bf Tors}}

%******************************************************************
\begin{abstract}
We introduce  a cohomology theory for  a class of projective
varieties over a finite field  coming  from  the canonical trace
on  a  $C^*$-algebra attached to  the variety.
Using the cohomology,  we prove  the rationality,  functional equation and the Betti numbers
conjectures for  the zeta function of the variety.

\vspace{7mm}


\vspace{5mm}

\end{abstract}

%**************************************************************************
\section{Introduction}
%***************************************************************************
The aim of our  note is a cohomology theory for projective varieties
over the field with $q=p^r$ elements.  Such a cohomology 
comes from the canonical trace on a $C^*$-algebra
attached  to  the  variety $V$;  this theory will be called a {\it trace
cohomology} and denoted  by $H^*_{tr}(V)$.   
The  trace cohomology  is much  parallel  to  the  $\ell$-adic cohomology  
$H^*_{et}(V; \Q_{\ell})$,  see  [Grothendieck 1968]  \cite{Gro1} and  [Hartshorne 1977,  pp.  453-457]  \cite{H}.
 Unlike  the $\ell$-adic cohomology,  it does not  depend on a prime  $\ell$
 and the  endomorphisms of  $H^*_{tr}(V)$ always have an integer trace.  
 While the $\ell$-adic  cohomology counts   isolated fixed points of the Frobenius endomorphism geometrically,   
 the trace cohomology does it  algebraically,  i.e.  taking into account the index $\pm 1$ of a  fixed point.  
 Moreover,  the eigenvalues of Frobenius endomorphism on $H^*_{et}(V; \Q_{\ell})$
are  complex  algebraic numbers of   the absolute value $q^{i\over 2}$,  yet   
such eigenvalues  are {\it real}  algebraic  on the trace cohomology $H^*_{tr}(V)$. 
 The  cohomology   groups $H^i_{tr}(V)$   are  truly  concrete and simple;  they  can
 be  found explicitly  in many important special cases,  e.g.  when  $V$ is an algebraic curve,
 see Section 4.   We shall pass to a detailed construction.

Denote by $V_{\Bbb C}$ an $n$-dimensional projective variety
over the field of complex numbers,
such that the reduction of $V_{\Bbb C}$  modulo the prime ideal
over $p$ is isomorphic to the variety $V:=V({\Bbb F}_q)$.
(In other words,   we assume that  $V$  can be  lifted to the characteristic zero  and 
the reduction modulo $p$ is functorial;    the corresponding  category
is described in   [Hartshorne 2010,  Theorem 22.1]  \cite{H2}.)
Let $B(V_{\Bbb C},  {\cal L},  \sigma)$ be the twisted homogeneous coordinate ring of
 projective variety $V_{\Bbb C}$,   where ${\cal L}$ is the invertible sheaf
of linear forms on $V_{\Bbb C}$ and $\sigma$ an automorphism of $V_{\Bbb C}$,   see
[Stafford \& van ~den ~Bergh 2001, p. 180]  \cite{StaVdb1} 
for the details.
The norm-closure of a self-adjoint representation of the  ring $B(V_{\Bbb C},  {\cal L},  \sigma)$
by the bounded linear operators on a Hilbert space ${\cal H}$  is a $C^*$-algebra,
see e.g. [Murphy 1990]   \cite{M}  for an introduction;   we 
call it a  {\it Serre $C^*$-algebra} of $V_{\Bbb C}$ and denote 
by ${\cal A}_V$.     Let ${\cal K}$ be the $C^*$-algebra of all compact
operators on ${\cal H}$.  We shall write $\tau: {\cal A}_V\otimes {\cal K}\to {\Bbb R}$
to denote   the canonical  normalized trace on  ${\cal A}_V\otimes {\cal K}$,   i.e. a positive linear functional
of norm $1$  such that $\tau(yx)=\tau(xy)$ for all $x,y\in {\cal A}_V\otimes {\cal K}$,  see  
 [Blackadar 1986, p. 31] \cite{B}.  
Because ${\cal A}_V$  is a crossed product $C^*$-algebra of the form
${\cal A}_V\cong C(V_{\Bbb C})\rtimes {\Bbb Z}$,   one can use  the Pimsner-Voiculescu 
six term exact sequence for the crossed products,  see  e.g.  [Blackadar 1986,  p. 83]  \cite{B}  for
 the details.  Thus   one gets the  short exact sequence of the algebraic $K$-groups:  
$0\to K_0(C(V_{\Bbb C}))\buildrel  i_*\over\to  K_0({\cal A}_V)\to K_1(C(V_{\Bbb C}))\to 0$, 
where   map  $i_*$  is induced by the natural  embedding of $C(V_{\Bbb C})$ 
into ${\cal A}_V$.   We  have $K_0(C(V_{\Bbb C}))\cong K^0(V_{\Bbb C})$ and 
$K_1(C(V_{\Bbb C}))\cong K^{-1}(V_{\Bbb C})$,  where $K^0$ and $K^{-1}$  are  the topological
$K$-groups of  $V_{\Bbb C}$, see  [Blackadar 1986, p. 80]  \cite{B}. 
By  the Chern character formula,  one gets
$K^0(V_{\Bbb C})\otimes {\Bbb Q}\cong H^{even}(V_{\Bbb C}; {\Bbb Q})$ and 
$K^{-1}(V_{\Bbb C})\otimes {\Bbb Q}\cong H^{odd}(V_{\Bbb C}; {\Bbb Q})$, 
where $H^{even}$  ($H^{odd}$)  is the direct sum of even (odd, resp.) 
cohomology groups of $V_{\Bbb C}$.  
Notice that $K_0({\cal A}_V\otimes {\cal K})\cong K_0({\cal A}_V)$ because
of  a stability of the $K_0$-group with respect to tensor products by the algebra 
${\cal K}$,  see e.g.   [Blackadar 1986,  p. 32]  \cite{B}.
One gets the   commutative diagram in Fig. 1, 
where $\tau_*$ denotes  a homomorphism  induced on $K_0$ by  the canonical  trace 
$\tau$ on the $C^*$-algebra  ${\cal A}_V\otimes {\cal K}$. 
Since  $H^{even}(V_{\Bbb C}):=\oplus_{i=0}^n H^{2i}(V_{\Bbb C})$ and  
$H^{odd}(V):=\oplus_{i=1}^n H^{2i-1}(V_{\Bbb C})$,   one gets  for each  $0\le i\le 2n$ 
 an injective  homomorphism   $ \tau_*:  ~H^i(V_{\Bbb C})\longrightarrow  {\Bbb R}$. 
%*******************************************************************
\begin{figure}[here]
%*******************************************************************
%*******************************************************************
\begin{picture}(300,100)(-40,0)
\put(160,72){\vector(0,-1){35}}
\put(80,65){\vector(2,-1){45}}
\put(240,65){\vector(-2,-1){45}}
\put(10,80){$ H^{even}(V_{\Bbb C})\otimes {\Bbb Q} 
\buildrel  i_*\over\longrightarrow  K_0({\cal A}_V\otimes{\cal K})\otimes {\Bbb Q} 
\longrightarrow H^{odd}(V_{\Bbb C})\otimes {\Bbb Q}$}
\put(167,55){$\tau_*$}
\put(157,20){${\Bbb R}$}
\end{picture}
%***********************************************************
%***********************************************************
\caption{The trace cohomology.}
\end{figure}
%*******************************************************************
 %********************************************************************************
\begin{dfn}\label{dfn1}
By an $i$-th trace cohomology  group $H^i_{tr}(V)$  of  variety  $V$   one 
understands the  abelian subgroup  of   ${\Bbb R}$ defined by the map $\tau_*$.
\end{dfn}
%*********************************************************************************
Notice that each  endomorphism  of  $H^i_{tr}(V)$ is  given   by 
a  real number  $\omega$,   such that $\omega H^i_{tr}(V)\subseteq H^i_{tr}(V)$; 
thus the ring $End~(H^i_{tr}(V))$ of all  endomorphisms of   $H^i_{tr}(V)$
is  commutative.   The $End~(H^i_{tr}(V))$ is a commutative subring of the ring  
$End~(H^i(V_{\Bbb C}))$  of all endomorphisms of the cohomology group $H^i(V_{\Bbb C})$.  
Moreover,  each regular map $f: V\to V$ corresponds to an algebraic map 
$f_{\Bbb C}: V_{\Bbb C}\to V_{\Bbb C}$  and, therefore,  to an endomorphism 
$\omega\in End~(H^i_{tr} (V))$. 
On the other hand,   it is easy to see that   $End~(H^i_{tr}(V))\cong {\Bbb Z}$ or
 $End~(H^i_{tr}(V))\otimes {\Bbb Q}$  is an algebraic number field.    
 In the latter case  $H^i_{tr}(V)\subset  End~(H^i_{tr}(V))\otimes {\Bbb Q}$,
see  [Manin 2004, Lemma 1.1.1]  \cite{Man1}   for the case of quadratic fields.  
We shall write  $tr~(\omega)$ to denote
the trace of an algebraic number $\omega\in End~(H^i_{tr} (V))$.  
 Our main results  are  as follows.
 %***************************************************************************
 \begin{thm}\label{thm0}
 The cardinality of  variety $V({\Bbb F}_{q})$  is given by the formula:
  %***************************************************************************
 \begin{equation}\label{eq1}
 |V({\Bbb F}_{q})|=1+q^n +\sum_{i=1}^{2n-1} (-1)^i ~tr~(\omega_i),
\end{equation}
%****************************************************************************
where $\omega_i\in End~(H^i_{tr}(V))$ is  generated  by the 
Frobenius map of $V({\Bbb F}_q)$.  
 \end{thm}
 %**************************************************************************
 %***********************************************************************
\begin{thm}\label{thm1}
The zeta function  $Z_V(t):=\exp\left(\sum_{r=1}^{\infty} {|V({\Bbb F}_{q^r})|\over r}t^r\right)$
 of $V({\Bbb F}_q)$  has the following properties:  

\medskip
(i)  $Z_V(t)={P_1(t)\dots P_{2n-1}(t)\over P_0(t)\dots P_{2n}(t)}$ is a rational function;

\smallskip
(ii) $Z_V(t)$ satisfies the functional equation  $Z_V\left({1\over q^nt}\right)=\pm q^{n{\chi(V_{\Bbb C})\over 2}} t^{\chi(V_{\Bbb C})} Z_V(t)$,
where $\chi(V_{\Bbb C})$ is the Euler-Poincar\'e characteristic of $V_{\Bbb C}$;

\smallskip
(iii)  $\deg~P_i(t)=\dim H^i(V_{\Bbb C})$. 
\end{thm}
%**************************************************************************
%**************************************************************************
\begin{rmk}
\textnormal{
Roughly speaking,  theorem \ref{thm1} says that the standard properties of the trace cohomology imply  all  Weil's  conjectures,
except  for    an  analog of the Riemann hypothesis   $|\alpha_{ij}|=q^{i\over 2}$ 
for  the roots  $\alpha_{ij}$  of  polynomials  $P_i(t)$,   see [Weil  1949,  p. 507]  \cite{Wei1};
the latter property   is proved in   \cite{Nik2}.
}
\end{rmk}
%********************************************************************************
The article is organized as follows.  Section 2 contains some useful
definitions and notation.  Theorems  \ref{thm0} and \ref{thm1} are  proved in Section 3.  
We calculate the trace cohomology for the   algebraic  (and elliptic, in particular)  curves  in Section 4.

%**************************************************************************
\section{Preliminaries}
%***************************************************************************
In this section  we briefly review the twisted homogeneous coordinate 
rings and the Serre $C^*$-algebras associated to projective varieties,
see   [Artin \& van den Bergh  1990]  \cite{ArtVdb1}) and 
 [Stafford \& van ~den ~Bergh 2001]  \cite{StaVdb1} for a detailed account. 
The $C^*$-algebras and their $K$-theory are covered in [Murphy 1990]  \cite{M}
and [Blackadar 1986]  \cite{B}, respectively.  The Serre $C^*$-algebras were 
introduced in \cite{Nik1}.

%**************************************************************************
\subsection{Twisted homogeneous coordinate rings}
%***************************************************************************
Let $V$ be a projective scheme over a field $k$, and let ${\cal L}$ 
be the invertible sheaf ${\cal O}_V(1)$ of linear forms on $V$.  Recall,
that the homogeneous coordinate ring of $V$ is a graded $k$-algebra, 
which is isomorphic to the algebra
%*************************************************************************
\begin{equation}\label{eq5}
B(V, {\cal L})=\bigoplus_{n\ge 0} H^0(V, ~{\cal L}^{\otimes n}). 
\end{equation}
%************************************************************************* 
Denote by $\Coh$ the category of quasi-coherent sheaves on a scheme $V$
and by $\Mod$ the category of graded left modules over a graded ring $B$.  
If $M=\oplus M_n$ and $M_n=0$ for $n>>0$, then the graded module
$M$ is called  right bounded.  The  direct limit  $M=\lim M_{\alpha}$
is called a torsion, if each $M_{\alpha}$ is a right bounded graded
module. Denote by $\Tors$ the full subcategory of $\Mod$ of the torsion
modules.  The following result is basic about the graded ring $B=B(V, {\cal L})$.   
%************************************************************************
\begin{lem}\label{lm1}
{\bf ([Serre 1955]  \cite{Ser1})}
\quad $\Mod~(B) ~/~\Tors \cong \Coh~(V).$
\end{lem}
%***************************************************************************
Let $\sigma$ be an automorphism of $V$.  The pullback of sheaf ${\cal L}$ 
along $\sigma$ will be denoted by ${\cal L}^{\sigma}$,  i.e. 
${\cal L}^{\sigma}(U):= {\cal L}(\sigma U)$ for every $U\subset V$. 
The graded $k$-algebra
 %*************************************************************************
\begin{equation}\label{eq7}
B(V, {\cal L}, \sigma)=\bigoplus_{n\ge 0} H^0(V, ~{\cal L}\otimes {\cal L}^{\sigma}\otimes\dots
\otimes  {\cal L}^{\sigma^{ n-1}}). 
\end{equation}
%************************************************************************* 
 is called a {\it twisted homogeneous coordinate ring} of scheme $V$;  
 notice that such a ring is non-commutative,  unless $\sigma$ is the trivial 
 automorphism.   The multiplication of sections is defined by the rule
 %***********************************************************************
 %\begin{equation}\label{eq8}
 $ab=a\otimes b^{\sigma^m}$,
 %\end{equation}
 %*************************************************************************
 whenever $a\in B_m$ and $b\in B_n$. 
 Given a pair $(V,\sigma)$ consisting of a Noetherian scheme $V$ and 
 an automorphism $\sigma$ of $V$,  an invertible sheaf ${\cal L}$ on $V$
 is called  $\sigma$-ample,  if for every coherent sheaf ${\cal F}$ on $V$,
 the cohomology group $H^q(V, ~{\cal L}\otimes {\cal L}^{\sigma}\otimes\dots
\otimes  {\cal L}^{\sigma^{ n-1}}\otimes {\cal F})$  vanishes for $q>0$ and
$n>>0$.  Notice,  that if $\sigma$ is trivial,  this definition is equivalent to the
usual definition of ample invertible sheaf [Serre 1955]  \cite{Ser1}.    
A  non-commutative generalization of the Serre theorem is as follows.
 %************************************************************************
\begin{lem}\label{lm2}
{\bf ([Artin \& van den Bergh  1990]  \cite{ArtVdb1})}
Let $\sigma: V\to V$ be an automorphism of a projective scheme $V$
over $k$  and let ${\cal L}$ be a $\sigma$-ample invertible sheaf on $V$. If
$B(V, {\cal L}, \sigma)$  is  the ring (\ref{eq7}),    then
%************************************************************************************
\begin{equation}
\Mod~(B(V, {\cal L}, \sigma)) ~/~\Tors \cong \Coh~(V).  
\end{equation}
%******************************************************************************
\end{lem}
%***************************************************************************

%**************************************************************************
\subsection{Serre $C^*$-algebras}
%***************************************************************************
Let $V$ be a projective scheme and $B(V, {\cal L}, \sigma)$ its
twisted homogeneous coordinate ring.   Let $R$ be a commutative 
graded ring, such that $V=Spec~(R)$.  Denote by $R[t,t^{-1}; \sigma]$
the ring of skew Laurent polynomials defined by the commutation relation
%************************************************************************************
%\begin{equation}\label{eq10}
$b^{\sigma}t=tb$
%\end{equation}
%******************************************************************************
for all $b\in R$, where $b^{\sigma}$ is the image of $b$ under automorphism 
$\sigma: V\to V$.  
 %************************************************************************
\begin{lem}\label{lm3}
{\bf ([Artin \& van den Bergh  1990]  \cite{ArtVdb1})}
$R[t,t^{-1}; \sigma]\cong B(V, {\cal L}, \sigma)$. 
\end{lem}
%***************************************************************************
Let ${\cal H}$ be a Hilbert space and   ${\cal B}({\cal H})$ the algebra of 
all  bounded linear  operators on  ${\cal H}$.
For a  ring of skew Laurent polynomials $R[t, t^{-1};  \sigma]$,  
we shall consider a homomorphism 
%*************************************************************************
\begin{equation}\label{eq2bis}
\rho: R[t, t^{-1};  \sigma]\longrightarrow {\cal B}({\cal H}). 
\end{equation}
%*************************************************************************
Recall  that algebra ${\cal B}({\cal H})$ is endowed  with a $\ast$-involution;
the involution comes from the scalar product on the Hilbert space ${\cal H}$. 
We shall call representation (\ref{eq2bis})  $\ast$-coherent,   if
(i)  $\rho(t)$ and $\rho(t^{-1})$ are unitary operators,  such that
$\rho^*(t)=\rho(t^{-1})$ and 
(ii) for all $b\in R$ it holds $(\rho^*(b))^{\sigma(\rho)}=\rho^*(b^{\sigma})$, 
where $\sigma(\rho)$ is an automorphism of  $\rho(R)$  induced by $\sigma$. 
Whenever  $B=R[t, t^{-1};  \sigma]$  admits a $\ast$-coherent representation,
$\rho(B)$ is a $\ast$-algebra;  the norm-closure of  $\rho(B)$  yields
a   $C^*$-algebra,  see e.g.  [Murphy 1990,  Section 2.1]  \cite{M}.  
 We shall refer to such  as   the  {\it Serre $C^*$-algebra} and denote it 
by  ${\cal A}_V$.

 Recall that if   ${\cal A}$ is  a $C^*$-algebra and $\sigma: G\to Aut~({\cal A})$ is a 
 continuous homomorphism of the locally compact group $G$ 
group,  then  the triple $({\cal A}, G, \sigma)$   defines a  $C^*$-algebra  called  a  crossed product and  denoted by 
${\cal A}\rtimes_{\sigma}G$;   we refer the reader to   [Williams  2007, pp. 47-54]  \cite{W} 
for the details.  It is not hard to see,   that ${\cal A}_V$  is a crossed product 
$C^*$-algebra of the form    ${\cal A}_V\cong C(V)\rtimes_{\sigma} {\Bbb Z}$,  where $C(V)$ is the
$C^*$-algebra of all continuous complex-valued functions on $V$
and $\sigma$  is a $\ast$-coherent  automorphism of  $V$.

%**************************************************************************
\section{Proofs}
%***************************************************************************
%**************************************************************************
\subsection{Proof of theorem \ref{thm0}}
%***************************************************************************
We shall prove a stronger result contained in the following lemma.
%************************************************************************
\begin{lem}\label{lm5}
The Lefschetz number of the Frobenius map $f_{\Bbb C}: V_{\Bbb C}\to V_{\Bbb C}$
is given by the formula:
%*******************************************************************************
\begin{equation}\label{eq6}
L(f_{\Bbb C})=1-q^n+\sum_{i=1}^{2n-1} (-1)^i ~tr~(\omega_i). 
\end{equation}
%********************************************************************************
\end{lem}
%***************************************************************************
{\it Proof.} 
Recall that the Lefschetz number of a continuous map $g_{\Bbb C}: V_{\Bbb C}\to V_{\Bbb C}$
is defined as 
%*******************************************************************************
\begin{equation}
L(g_{\Bbb C})=\sum_{i=0}^{2n} (-1)^i ~tr~(g_i^*), 
\end{equation}
%******************************************************************************** 
where $g_i^*: H^i(V_{\Bbb C})\to H^i(V_{\Bbb C})$ is an induced linear map  
of the cohomology.  Because $f^*_i$ is nothing but the matrix form of an endomorphism
$\omega_i\in End~(H^i_{tr}(V))$,  one gets 
%*******************************************************************************
\begin{equation}\label{eq8}
L(f_{\Bbb C})=\sum_{i=0}^{2n} (-1)^i ~tr~(\omega_i). 
\end{equation}
%********************************************************************************
We shall write equation (\ref{eq8}) in the form
%*******************************************************************************
\begin{equation}\label{eq9}
L(f_{\Bbb C})=tr~(\omega_0)+tr~(\omega_{2n})+\sum_{i=1}^{2n-1} (-1)^i ~tr~(\omega_i). 
\end{equation}
%********************************************************************************
It is known,  that $H^0(V_{\Bbb C})\cong {\Bbb Z}$ and $\omega_0=1$ is the trivial endomorphism;
thus $tr~(\omega_0)=1$.   Likewise,  $H^{2n}\cong {\Bbb Z}$,   but 
%*******************************************************************************
\begin{equation}\label{eq10}
\omega_{2n}=sgn ~[N(\omega_1)] ~q^n,
\end{equation}
%********************************************************************************
where $N(\bullet)$ is the norm of an algebraic number.  
It is known,  that the endomorphism  $\omega_1\in End~(H^1_{tr}(V))$
has the following matrix form
%*******************************************************************************
\begin{equation}\label{eq11}
q^{1\over 2}\left(\matrix{A & I\cr I & 0}\right),
\end{equation}
%********************************************************************************
where $A$ is a positive symmetric  and $I$ is the identity matrix,  see  (\cite{Nik1},
Lemma 3).  Thus
%*************************************************************
\begin{eqnarray}\label{eq12}
sgn~[N(\omega_1)] &= sgn ~\det  \left(\matrix{A & I\cr I & 0}\right)   =& \nonumber \\
 = sgn [-\det~(I^2)] &= -sgn \det~(I)=-1. &
\end{eqnarray}
%****************************************************************
Therefore,  from (\ref{eq10}) one obtains $\omega_{2n}=-q^n$;  in other words,   the Frobenius endomorphism
acts on $H^{2n}_{tr}(V)\cong {\Bbb Z}$  by multiplication on the negative integer $-q^n$. Clearly,
$tr~(\omega_{2n})=-q^n$ and the substitution of  these data  in (\ref{eq9}) gives us
%*******************************************************************************
\begin{equation}\label{eq13}
L(f_{\Bbb C})=1-q^n+\sum_{i=1}^{2n-1} (-1)^i ~tr~(\omega_i). 
\end{equation}
%******************************************************************************** 
Lemma \ref{lm5} follows.
$\square$

%***********************************************************************************
\begin{cor}\label{cr1}
The total number of  the index $-1$  fixed points of the Frobenius map  $f_{\Bbb C}$
is equal to $q^n$.  
\end{cor}
%***********************************************************************************
{\it Proof.}
It is known,  that 
 %***************************************************************************
 \begin{equation}\label{eqn14}
 |V({\Bbb F}_{q})|=1+q^n +\sum_{i=1}^{2n-1} (-1)^i ~tr~(Fr_i^*),
\end{equation}
%****************************************************************************
where $Fr_i^*:  H^i_{et}(V; \Q_{\ell})\to H^i_{et}(V; \Q_{\ell})$   is a linear
map on the $i$-th $\ell$-adic cohomology induced by the Frobenius endomorphism
of $V$,  see  [Hartshorne 1977, pp.  453-457]  \cite{H}.   But according to \cite{Nik1},   it holds
 %***************************************************************************
 \begin{equation}\label{eqn15}
 tr~(Fr_i^*)=tr~(\omega_i),\quad 1\le i\le 2n-1.  
\end{equation}
%****************************************************************************
Since $|Fix~(f_{\Bbb C})|=|V({\Bbb F}_q)|$,  one concludes from lemma \ref{lm5}
that the algebraic count $L(f_{\Bbb C})$ of the fixed points of $f_{\Bbb C}$ differs from its geometric
count $|Fix~(f_{\Bbb C})|$  by exactly $q^n$ points  of the index $-1$.   Corollary \ref{cr1}
is proved.
$\square$

\bigskip
Theorem \ref{thm0} follows formally from the equations (\ref{eqn14})  and (\ref{eqn15}).
$\square$

%**************************************************************************
\subsection{Proof of theorem \ref{thm1}}
%***************************************************************************
For the sake of clarity,   let us outline the main idea.  
Since  the trace cohomology accounts for the fixed points of the Frobenius 
map $f_{\Bbb C}$ algebraically (see corollary \ref{cr1}),   we shall deal with the corresponding Lefschetz 
zeta function
%***************************************************************************
 \begin{equation}\label{eqn16}
 Z_V^L(t):=\exp\left(\sum_{r=1}^{\infty} {L(f_{\Bbb C}^r)\over r}t^r\right)
\end{equation}
%****************************************************************************
and prove items (i)-(iii) for the $Z_V^L(t)$.  Because $f_{\Bbb C}: V_{\Bbb C}\to V_{\Bbb C}$
is the Anosov-type map,  one can use Smale's formulas linking $Z_V^L(t)$ and  $Z_V(t)$,
see [Smale 1967,  Proposition 4.14]  \cite{Sma1};  it will follow  that items  (i)-(iii)  are true for 
the function $Z_V(t)$ as well.  We shall pass to a detailed argument;  the following general  
lemma will be helpful.
%************************************************************************
\begin{lem}\label{lem5}
If  $f: V\to V$ is a regular map, then all eigenvalues $\lambda_{ij}$ of 
the corresponding endomorphisms  $\omega_i\in End~(H^i_{tr}(V))$
of the trace cohomology are real algebraic numbers.
\end{lem}
%***************************************************************************
{\it Proof.} 
Since the  endomorphisms of $H^i_{tr}(V)$ commute with each other, 
there exists  a basis of $H^i_{tr}(V)$,  such that each  endomorphism 
is given in this basis  by a symmetric integer matrix \cite{Nik1}.   But the spectrum of a real symmetric matrix is
known to  be totally real and  the eigenvalues of an integer matrix are algebraic numbers.
Lemma \ref{lem5} follows. 
$\square$

\bigskip
(i)  Let us prove  rationality of the function $Z_V^L(t)$ given by formula (\ref{eqn16}).  
Using lemma \ref{lm5},  one gets
%*************************************************************
\begin{eqnarray}\label{eq14}
\log~Z_V^L(t)  &=  ~\sum_{r=1}^{\infty}\left[1+(-q^n)^r+\sum_{i=1}^{2n-1} (-1)^i tr~(\omega_i^r)\right] {t^r\over r} =& \nonumber \\
 = \sum_{r=1}^{\infty} {t^r\over r} & +\sum_{r=1}^{\infty} {(-q^n t)^r\over r} +
 \sum_{r=1}^{\infty} \left(\sum_{i=1}^{2n-1} (-1)^i tr~(\omega_i^r)\right) {t^r\over r}. &
\end{eqnarray}
%****************************************************************
Taking into account the  well-known summation formulas  
$\sum_{r=1}^{\infty} {t^r\over r}=-\log(1-t)$ and  
$\sum_{r=1}^{\infty} {(-q^n t)^r\over r}=-\log(1+q^nt)$,  one can bring equation (\ref{eq14}) to
 the form
%********************************************************************************
\begin{equation}\label{eq15}
\log~Z_V^L(t)=-\log(1-t)(1+q^nt) +\sum_{i=1}^{2n-1} (-1)^i \sum_{r=1}^{\infty} tr~(\omega_i^r) {t^r\over r}.   
\end{equation}
%*******************************************************************************
On the other hand, it easy to see that 
%********************************************************************************
\begin{equation}\label{eq16}
tr~(\omega_i^r)=\lambda_1^r+\dots+\lambda_{b_i}^r,
\end{equation}
%*******************************************************************************
where $\lambda_j$  are the eigenvalues of the Frobenius endomorphism $\omega_i\in  End~(H^i_{tr}(V))$  
and   $b_i$ is the  $i$-th Betti number of $V_{\Bbb C}$.  
Thus one can bring  (\ref{eq15}) to the form  
%********************************************************************************
\begin{eqnarray}\label{eq17}
\log~Z_V^L(t) &= -\log(1-t)(1+q^nt) +&\nonumber\\
+\sum_{i=1}^{2n-1}  & (-1)^i \sum_{r=1}^{\infty} 
\left[{(\lambda_1 t)^r\over r}+\dots+{(\lambda_{b_i} t)^r\over r}\right]. & 
\end{eqnarray}
%*******************************************************************************
Using the summation formula  $\sum_{r=1}^{\infty} {(\lambda_j t)^r\over r}=-\log(1-\lambda_j t)$,
one gets from (\ref{eq17})
%********************************************************************************
\begin{eqnarray}\label{eq18}
\log~Z_V^L(t) &= -\log(1-t)(1+q^nt) +&\nonumber\\
 \sum_{i=1}^{2n-1} &   (-1)^{i+1}  \log\left[(1-\lambda_1t)\dots (1-\lambda_{b_i}t)\right]. &
\end{eqnarray}
%*******************************************************************************
Notice that the product $(1-\lambda_1t)\dots (1-\lambda_{b_i}t)$ is nothing 
but the characteristic polynomial $P_i(t)$ of the Frobenius endomorphism
on the trace cohomology $H^i_{tr}(V)$;   thus   one can write (\ref{eq18}) in the form
%********************************************************************************
\begin{equation}
\log~Z_V^L(t)=\log {P_1(t)\dots P_{2n-1}(t)\over (1-t) P_2(t)\dots P_{2n-2}(t) (1+q^nt)}.
\end{equation}
%*******************************************************************************
 Taking exponents in the last equation,  one obtains
%********************************************************************************
\begin{equation}\label{eqn23}
Z_V^L(t)={P_1(t)\dots P_{2n-1}(t)\over P_0(t)\dots P_{2n}(t)},
\end{equation}
%******************************************************************************* 
 where $P_0(t)=1-t$ and $P_{2n}(t)=1+q^nt$.    Thus $Z_V^L(t)$ is a rational function.

To prove rationality of $Z_V(t)$,  recall that a map $f_{\Bbb C}: V_{\Bbb C}\to V_{\Bbb C}$
is called {\it Anosov-type},   if there exist a (possibly singular) pair of orthogonal foliations
${\cal F}_u$ and ${\cal F}_s$ of $V_{\Bbb C}$ preserved by $f_{\Bbb C}$.   
(Note that our definition is more general than the standard and includes all  continuous 
maps $f_{\Bbb C}$.)   Consider the trace cohomology $H^1_{tr}(V)$ endowed with the Frobenius 
endomorphism $\omega_1\in End~(H^1_{tr}(V))$.  Let ${\cal F}_s$ be a foliation of $V_{\Bbb C}$,
whose holonomy (Plante) group is isomorphic to $H^1_{tr}(V)$.  
  Because $\omega_1 H^1_{tr}(V)\subset H^1_{tr}(V)$,  one concludes that ${\cal F}_s$ is an 
  invariant stable foliation of the map $f_{\Bbb C}$.  The unstable foliation ${\cal F}_u$ 
  can be constructed likewise.  
  Thus $f_{\Bbb C}$  is the Anosov-type map of the manifold $V_{\Bbb C}$.
  One can apply now (an extension of)  [Smale 1967, Proposition 4.14]  \cite{Sma1},   which  says  that 
one of the following formulas  must hold:
%*******************************************************************
\begin{equation}\label{eqn24}
\left\{
\begin{array}{lll}
Z_V(t) &=&  {1\over Z_V^L(t)},\\
&&\\
Z_V(t) &=&  Z_V^L(-t),\\
&&\\
Z_V(t) &=&  {1\over Z_V^L(-t)}.   
\end{array}
\right.
\end{equation}
%*****************************************************************
Since $Z_V^L(t)$ is known to be a rational function (\ref{eqn23}),  it follows from  Smale's formulas (\ref{eqn24}) 
that  $Z_V(t)$ is rational as well.  Item (i)  is proved.

\bigskip
(ii)  Recall that the cohomology $H^*(V_{\Bbb C})$ satisfies the Poincar\'e dulality;  the 
duality can be given by a pairing
%********************************************************************************
\begin{equation}\label{eq21}
H^i(V_{\Bbb C})\times H^{2n-i}(V_{\Bbb C})\longrightarrow H^{2n}(V_{\Bbb C})
\end{equation}
%******************************************************************************* 
obtained  from  the cup-product on $H^*(V_{\Bbb C})$.

Let $f:V\to V$ be the Frobenius endomorphism and $f_{\Bbb C}: V_{\Bbb C}\to V_{\Bbb C}$
the corresponding algebraic map of $V_{\Bbb C}$  and consider the action $(f^{2n}_{\Bbb C})^*$
on the pairing $\langle\bullet,\bullet\rangle$ given by (\ref{eq21}).  
Since $H^{2n}(V_{\Bbb C})\cong {\Bbb Z}$ and the linear map $(f^{2n}_{\Bbb C})^*$
multiplies $H^{2n}(V_{\Bbb C})$ by the constant $q^n$,  one gets 
%********************************************************************************
\begin{equation}\label{eq22}
\langle (f^{i}_{\Bbb C})^*x,  (f^{2n-i}_{\Bbb C})^*y\rangle=q^n\langle x,y\rangle,   
\end{equation}
%******************************************************************************* 
for all $x\in H^i(V_{\Bbb C})$ and all $y\in H^{2n-i}(V_{\Bbb C})$. 
Recall the linear algebra identities,  given e.g. in    [Hartshorne 1977, Lemma 4.3, p. 456]  \cite{H};   
then (\ref{eq22})  implies   the following formulas
%*******************************************************************
\begin{equation}\label{eq23}
\left\{
\begin{array}{rcl}
\det~(I-(f^i_{\Bbb C})^*t) &=& {(-1)^{b_i} (q^n)^{b_i} t^{b_i}\over \det~(f^{2n-i}_{\Bbb C})^*}
\det~\left[ I-{1\over q^nt}(f^{2n-i}_{\Bbb C})^*\right]\\
\det~(f^i_{\Bbb C})^* &=& {(q^n)^{b_i}\over \det~(f^{2n-i}_{\Bbb C})^*},
\end{array}
\right.
\end{equation}
%*****************************************************************
where $b_i=\dim~H^i(V_{\Bbb C})$ are the $i$-th Betti numbers.   
But $\det~(I-(f^i_{\Bbb C})^*t):=P_i(t)$ and 
$\det~\left[ I-{1\over q^nt}(f^{2n-i}_{\Bbb C})^*\right]:=P_{2n-i}\left({1\over q^nt}\right)$;
therefore, the first equation of (\ref{eq23})  yields   the identity
%********************************************************************************
\begin{equation}\label{eq24}
P_i(t)=  {(-1)^{b_i} (q^n)^{b_i} \over \det~(f^{2n-i}_{\Bbb C})^*} t^{b_i}P_{2n-i}\left({1\over q^nt}\right).
\end{equation}
%******************************************************************************* 
Let us calculate $Z_V^L\left({1\over q^nt}\right)$ using (\ref{eq24});   one gets the following
expression
%********************************************************************************
\begin{eqnarray}\label{eq25}
Z_V^L\left({1\over q^nt}\right) =& {P_1\left({1\over q^nt}\right)\dots P_{2n-1}\left({1\over q^nt}\right)
\over P_0\left({1\over q^nt}\right)\dots P_{2n}\left({1\over q^nt}\right)} =   &\nonumber\\
=   {P_1(t)\dots P_{2n-1}(t)\over P_0(t)\dots P_{2n}(t)}   t^{(b_0-b_1+\dots)} & (-1)^{(b_0-b_1+\dots)}
  {\det~(f^0_{\Bbb C})^*\dots \det~(f^{2n}_{\Bbb C})^*\over \det~(f^1_{\Bbb C})^*\dots \det~(f^{2n-1}_{\Bbb C})^*}.&
\end{eqnarray}
%*******************************************************************************
Note that $b_0-b_1+\dots=\chi(V_{\Bbb C})$   is the Euler-Poincar\'e 
characteristic  of $V_{\Bbb C}$.   From the second equation of (\ref{eq23}) one  obtains 
the identity $det~(f^i_{\Bbb C})^*\det~(f^{2n-i}_{\Bbb C})^*=(q^n)^{b_i}$.   Thus (\ref{eq25})
can be written in the form 
%********************************************************************************
\begin{eqnarray}\label{eq26}
Z_V^L\left({1\over q^nt}\right)=&  t^{\chi(V_{\Bbb C})} (-1)^{\chi(V_{\Bbb C})} 
{(q^n)^{{1\over 2}(b_0+\dots b_{2n})}\over (q^n)^{{1\over 2}(b_1+\dots+b_{2n-1})}} Z_V^LV(t)=&\nonumber\\
= & t^{\chi(V_{\Bbb C})} (-1)^{\chi(V_{\Bbb C})}  (q^n)^{{1\over 2}\chi(V_{\Bbb C})} Z_V^L(t).&
\end{eqnarray}
%******************************************************************************* 
Taking into account  $(-1)^{-\chi(V_{\Bbb C})}=\pm 1$,  one gets 
a  functional equation for $Z_V^L(t)$.   We encourage the reader to verify using formulas
(\ref{eqn24}),  that the same equation holds for the function $Z_V(t)$.  
 Item (ii) of theorem \ref{thm1} is proved.

\bigskip
(iii)   To prove the Betti numbers conjecture,  notice that equality (\ref{eq16}) implies that
$\deg P_i(t)=\dim H^i_{tr}(V)$.   But  $\dim ~H^i_{tr}(V)=\dim ~H^i(V_{\Bbb C})$ by the definition
of trace cohomoogy;  thus $\deg P_i(t)=\dim H^i (V_{\Bbb C})$ for the polynomials $P_i(t)$
in formula (\ref{eqn23}).  Again,  the reader can verify using (\ref{eqn24}),  that the same 
relationship holds for the polynomials representing the rational function $Z_V(t)$.    Item (iii) is proved.

\bigskip
This argument completes the proof of theorem \ref{thm1}
$\square$

%**************************************************************************
\section{Examples}
%***************************************************************************
The  groups $H^i_{tr}(V)$  are  truly   concrete and simple;
in this section we calculate   the trace cohomology for  $n=1$,  i.e.  when $V$
is a smooth  algebraic curve.  In particular,  we find the cardinality of  the set ${\cal E}({\Bbb F}_{q})$ 
 obtained by the reduction modulo $q$   of   an elliptic curve  with  complex multiplication. 
The reader can verify, that the lifting condition for $V$ is satisfied, see footnote 1.  
%***********************************************************************
\begin{exm}\label{exm1}
\textnormal{
The trace cohomology of smooth  algebraic curve ${\cal C}({\Bbb F}_q)$ of
genus $g\ge 1$ is given by the formulas: 
%*******************************************************************
\begin{equation}\label{eq32}
\left\{
\begin{array}{lll}
H_{tr}^0({\cal C}) &\cong&  {\Bbb Z},\\
&&\\
H_{tr}^1({\cal C}) &\cong&  {\Bbb Z}+{\Bbb Z}\theta_1+\dots+{\Bbb Z}\theta_{2g-1},\\
&&\\
H_{tr}^2({\cal C}) &\cong&  {\Bbb Z},   
\end{array}
\right.
\end{equation}
%*****************************************************************
where $\theta_i\in {\Bbb R}$ are  algebraically independent  integers
of a number field of degree $2g$.  
}
\end{exm}
%**************************************************************************
{\it Proof.}
It is known that the Serre $C^*$-algebra of the (generic) complex algebraic curve ${\cal C}$
is isomorphic to a  {\it toric} $AF$-algebra  ${\Bbb A}_{\theta}$,  
see \cite{Nik3} for the notation and details.   Moreover,  up to  a scaling constant $\mu>0$,
it holds
%******************************************************************************************
\begin{equation}\label{eq33}
\tau_*(K_0({\Bbb A}_{\theta}\otimes {\cal K}))=
\cases{{\Bbb Z}+{\Bbb Z}\theta_1 & \hbox{if} $g=1$\cr
             {\Bbb Z}+{\Bbb Z}\theta_1+\dots+{\Bbb Z}\theta_{6g-7} & \hbox{if} $g>1$,}
\end{equation}
%***************************************************************************************** 
where constants $\theta_i\in {\Bbb R}$ parametrize the moduli (Teichm\"uller) space of curve ${\cal C}$,
{\it ibid.}   If ${\cal C}$ is defined over a number field $k$,  then each $\theta_i$ is algebraic and 
their total number   is equal to  $2g-1$.  (Indeed,  since $Gal~(\bar k ~|~ k)$ acts on the torsion
points of ${\cal C}(k)$,  it is easy to see that the endomorphism ring of ${\cal C}(k)$ is non-trivial.  
Because such a ring  is isomorphic to the endomorphism ring of jacobian $Jac~{\cal C}$ 
and $\dim_{\Bbb C} Jac~{\cal C}=g$,  one concludes that $End ~{\cal C}(k)$ is a ${\Bbb Z}$-module
of rank $2g$ and each $\theta_i$ is an algebraic number.)  After scaling by a constant $\mu>0$,   one gets   
%******************************************************************************************
\begin{equation}\label{eq34}
H^1_{tr}({\cal C}):=\tau_*(K_0({\Bbb A}_{\theta}\otimes {\cal K}))=
{\Bbb Z}+{\Bbb Z}\theta_1+\dots+{\Bbb Z}\theta_{2g-1}
\end{equation}
%***************************************************************************************** 
Because $H^0({\cal C})\cong H^2({\cal C})\cong {\Bbb Z}$,
one obtains   the rest of formulas (\ref{eq32}).
$\square$

%*********************************************************************
\begin{rmk}\label{rmk10}
\textnormal{
Using theorem \ref{thm0},  one gets the formula
%********************************************************************************
\begin{equation}\label{eq30}
|{\cal C}({\Bbb F}_{q})|=  1+q-tr~(\omega)= 1+q-\sum_{i=1}^{2g}\lambda_i, 
\end{equation}
%******************************************************************************* 
where $\lambda_i$ are  real eigenvalues of the Frobenius endomorphism 
$\omega\in End~(H^1_{tr}({\cal C}))$.      Note that 
%********************************************************************************
\begin{equation}
\lambda_1+\dots+\lambda_{2g}=\alpha_1+\dots+\alpha_{2g},
\end{equation}
%******************************************************************************* 
where $\alpha_i$ are the eigenvalues of the Frobenius endomorphism of  $H^1_{et}({\cal C}; \Q_{\ell})$.
However,  there is no trace cohomology analog of the  classical  formula
%********************************************************************************
\begin{equation}
|{\cal C}({\Bbb F}_{q^r})|=   1+q^r-\sum_{i=1}^{2g}\alpha_i^r,
\end{equation}
%******************************************************************************* 
unless $r=1$;   this difference is due to an algebraic count of the fixed 
points by the trace cohomology. 
}
\end{rmk}
%*********************************************************************
%***********************************************************************
\begin{exm}\label{exm2}
\textnormal{
The case $g=1$ is particularly  instructive;  for the sake of clarity,  we shall consider
elliptic curves having complex multiplication.  Let ${\cal E}({\Bbb F}_q)$  
be the reduction modulo $q$ of an elliptic with complex multiplication 
by the ring of integers of an imaginary quadratic field ${\Bbb Q}(\sqrt{-d})$,
see e.g.  [Silverman 1994, Chapter 2]   \cite{S}.  It is known,  that in this case 
the trace cohomology formulas (\ref{eq32})  take the form  
%*******************************************************************
\begin{equation}\label{eq38}
\left\{
\begin{array}{lll}
H_{tr}^0({\cal E}({\Bbb F}_q)) &\cong&  {\Bbb Z},\\
&&\\
H_{tr}^1({\cal E}({\Bbb F}_q)) &\cong&  {\Bbb Z}+{\Bbb Z}\sqrt{d},\\
&&\\
H_{tr}^2({\cal E}({\Bbb F}_q)) &\cong&  {\Bbb Z}.     
\end{array}
\right.
\end{equation}
%*****************************************************************
We shall denote by  $\psi({\goth P})\in {\Bbb Q}(\sqrt{-d})$  the Gr\"ossencharacter
of  the prime ideal ${\goth P}$ over $p$,   see  [Silverman 1994, p. 174]   \cite{S}.
It is easy to see, that in this case  the Frobenius endomorphism $\omega\in End~(H_{tr}^1({\cal E}({\Bbb F}_q)))$
is given by the formula
%*******************************************************************
\begin{equation}\label{eq39}
\omega={1\over 2}\left[\psi({\goth P})+\overline{\psi({\goth P})}\right] + 
{1\over 2}\sqrt{\left(\psi({\goth P})+\overline{\psi({\goth P})}\right)^2+4q}   
\end{equation}
%*****************************************************************
and the corresponding eigenvalues 
%******************************************************************************************
\begin{equation}\label{eq40}
\left\{
\begin{array}{lll}
\lambda_1 &=& \omega= {1\over 2}\left[\psi({\goth P})+\overline{\psi({\goth P})}\right] + 
{1\over 2}\sqrt{\left(\psi({\goth P})+\overline{\psi({\goth P})}\right)^2+4q},  \\
\lambda_2 &=& \bar\omega= {1\over 2}\left[\psi({\goth P})+\overline{\psi({\goth P})}\right] -
{1\over 2}\sqrt{\left(\psi({\goth P})+\overline{\psi({\goth P})}\right)^2+4q}  .
\end{array}
\right.
\end{equation}
%*****************************************************************************************
Using formula (\ref{eq30}),  one  gets the following equation 
%********************************************************************************
\begin{equation}\label{eqn41}
|{\cal E}({\Bbb F}_{q})|=  1-(\lambda_1+\lambda_2)+q=1-\psi({\goth P})-\overline{\psi({\goth P})} +q,
\end{equation}
%******************************************************************************* 
which coincides with the well-known expression for  $|{\cal E}({\Bbb F}_{q})|$ in terms 
of the Gr\"ossencharacter,  see e.g.  [Silverman 1994, p. 175]   \cite{S}. 
}
\end{exm}

\vskip1cm

\textsc{Department of Mathematics and Computer Science, St.~John's University, 8000 Utopia Parkway,  
New York,  NY 11439, United States;} ~\textsc{E-mail:} {\sf igor.v.nikolaev@gmail.com}


\begin{thebibliography}{100}

\bibitem{ArtVdb1}
M.~Artin and M. ~van den Bergh,  Twisted homogeneous coordinate
rings, J. of Algebra 133 (1990), 249-271. 



\bibitem{B}
B.~Blackadar, $K$-Theory for Operator Algebras, MSRI Publications,
Springer, 1986



\bibitem{Gro1}
A.~Grothendieck,  Standard conjectures on algebraic cycles,
in:  Algebraic Geometry,  Internat. Colloq. Tata Inst. Fund. Res., Bombay, 1968.


\bibitem{H}
R.~Hartshorne, Algebraic Geometry, GTM 52, Springer, 1977.



\bibitem{H2}
R.~Hartshorne, Deformation Theory, GTM 257, Springer, 2010.


\bibitem{Man1}
Yu.~I.~Manin, Real multiplication and noncommutative geometry,
in ``Legacy of Niels Hendrik Abel'', 685-727, Springer, 2004.


\bibitem{M}
G.~J.~Murphy,  $C^*$-Algebras and Operator Theory, Academic Press, 1990. 



\bibitem{Nik3}
I.~Nikolaev,   Noncommutative geometry of algebraic curves, 
Proc. Amer. Math. Soc. 137 (2009), 3283-3290.


\bibitem{Nik1}
I.~Nikolaev,  On  traces of  Frobenius   endomorphisms,  Finite Fields 
and their Applications 25 (2014),  270-279.   


\bibitem{Nik2}
I.~Nikolaev,  On trace cohomology, {\sf arXiv:1407.3982}



\bibitem{Ser1}
J.~P.~Serre,  Fasceaux alg\'ebriques coh\'erents, Ann. of Math. 61 (1955), 
197-278.  


\bibitem{S}
J.~H.~Silverman, Advanced Topics in the Arithmetic of Elliptic Curves,
GTM 151, Springer 1994.



\bibitem{Sma1}
S.~Smale,  Differentiable dynamical systems, Bull. Amer. Math. Soc. 73 (1967),
747-817. 


\bibitem{StaVdb1}
J.~T.~Stafford and M.~van ~den ~Bergh, Noncommutative curves and noncommutative
surfaces, Bull. Amer. Math. Soc. 38 (2001), 171-216. 


\bibitem{Wei1}
A.~Weil,  Numbers of solutions of equations in finite fields,   Bull. Amer. Math. Soc. 55
(1949), 497-508.   


\bibitem{W}
D.~P.~Williams,  Crossed Products of $C^*$-Algebras, Math. Surveys and 
Monographs, Vol. 134,  Amer. Math. Soc.  2007.
   


\end{thebibliography}
\end{document}